\documentclass{svmult}

\usepackage[sort]{cite}       
\usepackage{amssymb}          
\usepackage{latexsym}         
\usepackage{amsmath}          
\usepackage{amscd}            
\usepackage[latin1]{inputenc} 
\usepackage{eucal}            
\usepackage{bbm}              
\usepackage{theorem}          
\usepackage{mathrsfs}         
\usepackage{stmaryrd}         

\renewcommand{\mathbb}[1]{\mathbbm{#1}} 

\newcounter{comment}

\newcommand{\image}      {\operatorname{{\mathrm{im}}}} 
\newcommand{\Lie}        {\operatorname{\mathscr{L}}}    
\newcommand{\cc}[1]      {\overline{{#1}}}              
\newcommand{\id}         {\operatorname{\mathsf{id}}}   
  
\newcommand{\Pol}        {\operatorname{\mathrm{Pol}}}   
\newcommand{\tr}         {\operatorname{\mathsf{tr}}}    
\newcommand{\ad}         {\operatorname{\mathrm{ad}}}

\newcommand{\End}        {\operatorname{\mathsf{End}}}   
\newcommand{\SP}[1]      {\left\langle{#1}\right\rangle}

\newcommand{\red}        {\mathrm{red}}
\newcommand{\opp}        {\mathrm{opp}}
\newcommand{\hor}        {\mathrm{hor}}                
\newcommand{\pr}         {\mathrm{pr}}

\newcommand{\const}      {\mathit{const}}

\newcommand{\Sym}        {\mathrm{S}}
\newcommand{\Schouten}[1]{\left\llbracket{#1}\right\rrbracket}
\newcommand{\Diffop}     {\operatorname{\mathrm{Diffop}}}
\newcommand{\SymD}       {\mathsf{D}}
\newcommand{\divergence} {\operatorname{\mathrm{div}}}
\newcommand{\ins}        {\operatorname{\mathrm{i}}}

\newcommand{\inss}       {\operatorname{\mathrm{i}_{\scriptscriptstyle\mathrm{s}}}}
\newcommand{\Ric}        {\operatorname{\mathrm{Ric}}} 
\newcommand{\DeltaRic}   {\operatorname{\Delta^{\mathrm{Ric}}}}
\newcommand{\nablaM}     {\nabla^M}

\newcommand{\starred}    {\mathbin{\star_\red}}
\newcommand{\weylrep}    {\operatorname{\varrho_{\scriptscriptstyle\mathrm{Weyl}}}}
\newcommand{\stdrep}     {\operatorname{\varrho_{\scriptscriptstyle\mathrm{Std}}}}
\newcommand{\krep}       {\operatorname{\varrho_{\kappa}}}

\newcommand{\qkoszul}    {\operatorname{\boldsymbol{\partial}}}
\newcommand{\qI}        {\boldsymbol{\mathcal{I}}_\Sigma}
\newcommand{\qB}        {\boldsymbol{\mathcal{B}}_\Sigma}

\newcommand{\qiota}      {\boldsymbol{\iota^*}}
\newcommand{\qh}         {\boldsymbol{h}}

\begin{document}

%
%

\title*{Symplectic Connections of Ricci Type and Star Products}

\author{Michel Cahen\inst{1}
  \and
  Simone Gutt\inst{1}\inst{2}
  \and
  Stefan Waldmann\inst{3}
}
 
\institute{D{\'e}partement de Math{\'e}matique \\
  Universit{\'e} Libre de Bruxelles \\
  Campus Plaine, C. P. 218 \\
  Boulevard du Triomphe \\
  B-1050 Bruxelles \\
  Belgium\\
  \texttt{mcahen@ulb.ac.be}, \texttt{sgutt@ulb.ac.be}
  \and 
  Universit\'e Paul Verlaine - Metz\\
  D{\'e}partement de Math{\'e}matique LMAM\\
  Ile du Saulcy \\
  F-57045~Metz Cedex 01\\
  France\\
  \texttt{gutt@poncelet.univ-metz.fr}
  \and
  Fakult{\"a}t f{\"u}r Mathematik und Physik \\
  Albert-Ludwigs-Universit{\"a}t Freiburg \\
  Physikalisches Institut \\
  Hermann Herder Strasse 3 \\
  D 79104 Freiburg \\
  Germany\\
  \texttt{Stefan.Waldmann@physik.uni-freiburg.de}
}
    
\date{September 2007}

\maketitle

\begin{abstract}
    In this article we relate the construction of  
    Ricci type symplectic connections by reduction 
    to the construction of star product by reduction 
  yielding rather explicit descriptions for the
    star product on the reduced space.
\end{abstract}

%
%

\section*{Introduction}
\label{sec:Introduction}

Deformation quantization \cite{bayen.et.al:1978a} is a formal
deformation ---in the sense of Murray Gerstenhaber
\cite{gerstenhaber:1964a}--- of the algebraic structure of the space
of smooth functions on a manifold $M$; it yields at first order in the
deformation parameter a Poisson structure on $M$. When this Poisson
structure is non-degenerate, i.e. when the manifold is symplectic, 
deformation quantization at second order yields a symplectic
connection on $M$ \cite{gutt.rawnsley:2003a}.

On a symplectic manifold $(M, \omega)$, symplectic connections always
exist but are not unique.  The curvature $R^\nabla$ of such a
connection $\nabla$ splits \cite{vaisman:1985a} under the action of
the symplectic group, (when the dimension of $M$ is at least $4$),
into two irreducible components $R^\nabla=E^\nabla+W^\nabla$ with
$E^\nabla$ completely determined by the Ricci tensor of the
connection. A symplectic connection is said to be of Ricci type if
$R^\nabla=E^\nabla$.

Marsden-Weinstein reduction, see e.g.
\cite[Sect.~4.3]{abraham.marsden:1985a}, is a method in symplectic
geometry to construct a symplectic manifold $(M,\omega)$
---called the reduced space--- from a bigger one $(P, \mu)$ and the
extra data of a coisotropic submanifold $C$ in $P$.

Under some further asumptions, one can sometimes reduce connections
\cite{baguis.cahen:2003a, baguis.cahen:2001a}, i.e. define a
symplectic connection $\nabla^M$ on $M$ from a connection $\nabla^P$
on $P$.  In particular, any Ricci-type connection on a 
simply connected $2n(\ge 4)$
dimensional manifold can be obtained
\cite{cahen.gutt.schwachhoefer:2005a} by reduction from a flat
connection $\nabla$ on a $2n+2$ dimensional manifold $P$, with the
coisotropic codimension $1$ submanifold $C$ defined by the zero set of
a function $F: P \rightarrow \mathbb{R}$ whose third covariant
derivative vanishes.

One way to describe the algebra of functions on the reduced space is
the use of BRST methods. Jim Stasheff participated actively to the
development of this point of view, see e.g.
\cite{fisch.henneaux.stasheff.teitelboim:1989a} among many others.
Reduction of deformation quantization has been studied by various
authors \cite{fedosov:1998a, bordemann.herbig.waldmann:2000a,
  bordemann:2004a:pre, cattaneo.felder:2005a:pre, gloessner:1998a:pre,
  cattaneo.felder:2004a}. In particular, a quantized version of BRST
methods was introduced in \cite{bordemann.herbig.waldmann:2000a} to
construct reduction in deformation quantization.  We use here those
methods to define a star product on any symplectic manifold endowed
with a Ricci type connection.

In this context we rely both on the work of Murray and on the work of
Jim and we are very happy and honoured to dedicate this to them.

{\small{In section \ref{sec:Riccitypeconnections}, we recall some
    basic properties of Ricci type symplectic connections. Section
    \ref{sec:RicciOperator} introduces a natural differential operator
    of order $2$ on a symplectic manifold endowed with a symplectic
    connection. In section \ref{sec:WeylMoyalStarProduct}, we recall
    the expression of the Weyl-Moyal star product on a symplectic
    manifold with a flat connection.  Section
    \ref{sec:ReductionStarProducts} explains the construction of
    reduced star product in our context where the reduced space
    $(M,\omega)$ is a symplectic manifold of dimension $2n\ge 4$,
    endowed with a Ricci-type connection, and where the big space
    $(P,\mu)$ is of dimension $2n+2$ with a flat connection and the
    related Weyl-Moyal star product.  In section \ref{sec:Properties},
    we show some properties of the reduced star product, in particular
    that the connection defined by the reduced star product is the
    Ricci type connection.}}

%
%

\section{Preliminary Results on Ricci Type Connections}\label{sec:Riccitypeconnections}
\label{sec:RicciTypeConnections}

In this section we recall some basic properties of Ricci type
symplectic connections to explain our notation. We follow essentially
\cite{ cahen.gutt.schwachhoefer:2005a,
  cahen.schwachhoefer:2004a} and refer for further details
  to the expository paper  \cite{bieliavsky.et.al:2006a}.

Let $(M, \omega)$ be a $2n\ge 4$-dimensional symplectic manifold which
allows for a Ricci type symplectic connection: this is a symplectic
 connection $\nabla^M$ (i.e. a linear torsion-free connection so that
 the symplectic $2$-form $\omega$ is parallel) such that the curvature tensor $R$ is
entirely determined by its Ricci tensor $\Ric$. Precisely, for
$X, Y \in \Gamma^\infty(TM)$ we have
\begin{align}
    \nonumber
    R(X, Y) = - \frac{1}{2n+1} 
    &\Big(
    -2 \omega(X, Y) \varrho 
    - \varrho(Y) \otimes X^\flat \\
    &+ \varrho(X) \otimes Y^\flat
    - X \otimes (\varrho(Y))^\flat 
    + Y \otimes (\varrho(X))^\flat
    \Big),
    \label{eq:RforRicciType}
\end{align}
where $X^\flat = \ins_X \omega$ as usual and $\varrho$ is the Ricci
endomorphism defined by
\begin{equation}
    \label{eq:varrhoDef}
    \Ric(X, Y) = \omega(X, \varrho(Y)).
\end{equation}
It follows that there exists a vector field $\mathsf{U} \in
\Gamma^\infty(TM)$, a function $\mathsf{f} \in C^\infty(M)$ and a
constant $\mathsf{K}$ such that the following identities hold
\begin{equation}
    \label{eq:nablaXvarrho}
    \nabla^M_X \varrho = 
    - \frac{1}{2n+1} \left(
        X \otimes \mathsf{U}^\flat + \mathsf{U} \otimes X^\flat
        \right),
\end{equation}
\begin{equation}
    \label{eq:nablaXU}
    \nabla^M_X \mathsf{U} =
    - \frac{2n+1}{2(n+1)} \varrho^2X + \mathsf{f}X,
\end{equation}
\begin{equation}
    \label{eq:tracevarrho}
    \tr(\varrho^2) + \frac{4(n+1)}{2n+1} \mathsf{f} = \mathsf{K}.
\end{equation}

One of the fundamental properties of such a Ricci type connection is
that $(M, \omega, \nabla^M)$ can be obtained from a Marsden-Weinstein
reduction out of a $(2n+2)$-dimensional symplectic manifold $(P, \mu)$
which is equipped with a \emph{flat} symplectic torsion-free
connection $\nabla$. In fact, we have to assume that $M$ is
simply-connected;  then by \cite{cahen.gutt.schwachhoefer:2005a} there exists
a $(2n+1)$-dimensional manifold $\Sigma$ with a surjective submersion
\begin{equation}
    \label{eq:piSigmaM}
    \pi: \Sigma \longrightarrow M
\end{equation}
together with a contact one-form $\alpha \in
\Gamma^\infty(T^*\Sigma)$, i.e. $\alpha \wedge (\D\alpha)^n$ is
nowhere vanishing  with $\pi^*\omega = \D \alpha$, whose Reeb vector
field $Z \in \Gamma^\infty(T\Sigma)$, defined by $\alpha(Z) = 1$ and
$\ins_Z\D\alpha = 0$, has a flow such that \eqref{eq:piSigmaM} is the
quotient onto the orbit space with respect to this flow. 
This manifold  $\Sigma$ is constructed as the holonomy bundle
over $M$ for a connection defined on an extension
of the frame bundle. It is a
$\mathbb{R}$ or $\mathbb{S}^1$ principal bundle over $M$ with
connection one-form $\alpha$.\\
Then we
consider $P = \Sigma \times \mathbb{R}$ with $\pr_1: P \longrightarrow
\Sigma$ being the canonical projection and $\iota: \Sigma
\longrightarrow P$ being the embedding of $\Sigma$ as $\Sigma \times
\{0\}$. The coordinate along $\mathbb{R}$ is denoted by $s$ and we set
$S = \frac{\partial}{\partial s} \in \Gamma^\infty(TP)$. On $P$ one
has the following exact symplectic form
\begin{equation}
    \label{eq:muDef}
    \mu = \D(\E^s \pr_1^*\alpha) 
    = \E^s \D s \wedge \pr_1^*\alpha + \E^s \D \pr_1^*\alpha.
\end{equation}
Thanks to the Cartesian product structure we can lift vector fields on
$\Sigma$ canonically to $P$. In particular, the lift $E$ of the Reeb
vector field $Z$ (defined by $ds(E)=0$ and $Tpr_1( E)=Z $) 
turns out to be Hamiltonian $E = - X_H$ with
\begin{equation}
    \label{eq:HDef}
    H = \E^s \in C^\infty(P).
\end{equation}
We have
\begin{equation}
    \label{eq:SXHSHSmu}
    \Lie_S X_H = 0,
    \quad
    \Lie_S H = H,
    \quad
    \textrm{and}
    \quad
    \Lie_S \mu = \mu,
\end{equation}
whence in particular $S$ is a conformally symplectic vector
field. Moreover, $\mu(X_H, S) = H$ and we can rewrite $\mu$ as
\begin{equation}
    \label{eq:muAgain}
    \mu = \D H \wedge \pr_1^*\alpha + H \pr^*\omega,
\end{equation}
where $\pr = \pi \circ \pr_1: P \longrightarrow M$. Then $(M, \omega)$
is the (Marsden-Weinstein) reduced phase space of $(P, \mu)$ with
respect to the Hamiltonian flow of $H$ at momentum value $H = 1$,
since indeed $\Sigma = H^{-1}(\{1\})$ and $\iota^*\mu = \pi^*\omega$
by \eqref{eq:muAgain}.

Using the contact form $\alpha$ we can lift vector fields $X \in
\Gamma^\infty(TM)$ horizontally to vector fields $\overline{X} \in
\Gamma^\infty(T\Sigma)$ by the condition
\begin{equation}
    \label{eq:FirstHorizontalLift}
    T\pi \circ \overline{X} =  X \circ \pi
    \quad
    \textrm{and}
    \quad
    \alpha(\overline{X}) = 0.
\end{equation}
 Since  $\D \alpha =
\pi^*\omega$ we have $[\overline{X}, \overline{Y}] = \overline{[X, Y]}
- \pi^*(\omega(X, Y)) Z$. Using also the canonical lift to $P$, we can
lift $X \in \Gamma^\infty(TM)$ horizontally to $X^\hor \in
\Gamma^\infty(TP)$, now subject to the conditions
\begin{equation}
    \label{eq:horizontalLift}
    T\pr \circ X^\hor = X \circ \pr
    \quad
    \textrm{and}
    \quad
    \pr_1^*\alpha(X^\hor) = 0 = \D s(X^\hor).
\end{equation}
We have 
\begin{equation}
    \label{eq:XhorYhorLiebracket}
    \left[X^\hor, Y^\hor\right]
    =
    [X, Y]^\hor + \pr^*(\omega(X, Y)) X_H
\end{equation}
as well as
\begin{equation}
    \label{eq:LieSXhorLieXHXhor}
    \left[S, X^\hor\right] = 0 = \left[X_H, X^\hor\right].
\end{equation}

We shall speak of ``invariance'' always with respect to the flow of
$X_H$ (or $Z$ on $\Sigma$, respectively) and of ``homogeneity'' always
with respect to the conformally symplectic vector field $S$, e.g. a
differential operator $D$ on $P$ is called homogeneous of degree $k
\in \mathbb{Z}$ if $[\Lie_S, D] = kD$, etc.

Denote the Poisson tensor on $M$ by $\Lambda_M$ and the one on $P$ by
$\Lambda_P$, respectively. We can also extend the horizontal lift
to bivectors as usual. Since we have curvature, $\Lambda_M^\hor$ is no
longer a Poisson bivector, instead one finds
\begin{equation}
    \label{eq:SchoutenLambdahor}
    \Schouten{\Lambda_M^\hor, \Lambda_M^\hor}
    =
    -2 \Lambda_M^\hor \wedge X_H.
\end{equation}
From the Definition (\ref{eq:muDef}) of the symplectic $2$-form $\mu$
on $P$, we have the relation
\begin{equation}
    \label{eq:PoissonP}
    \Lambda_P = 
    \frac{1}{H} \left(
        \Lambda_M^\hor + S \wedge X_H
    \right).
\end{equation}
In particular, for $u, v \in C^\infty(M)$ we find for the Poisson
brackets
\begin{equation}
    \label{eq:PoissonBracketuv}
    \{\pr^*u, \pr^*v\} _P= \frac{1}{H} \pr^*\{u, v\}_M.
\end{equation}

We are now in the position to define the flat connection $\nabla$ on
$P$ by specifying it on horizontal lifts, on $S$ and on $X_H$. One
defines \cite{cahen.gutt.schwachhoefer:2005a}
\begin{equation}
    \label{eq:nablaXhorYhor}
    \nabla_{X^\hor} Y^\hor 
    = \left(\nabla^M_X Y\right)^\hor 
    + \frac{1}{2} \pr^*(\omega(X, Y)) X_H + \pr^*(t(X, Y))S,
\end{equation}
\begin{equation}
    \label{eq:nablaXhorXH}
    \nabla_{X^\hor} X_H = \nabla_{X_H} X^\hor
    = \left(\tau(X)\right)^\hor
    - \pr^*(\omega(X, \mathsf{V})) S,
\end{equation}
\begin{equation}
    \label{eq:nablaXhorS}
    \nabla_{X^\hor} S = \nabla_S X^\hor = \frac{1}{2} X^\hor,
\end{equation}
\begin{equation}
    \label{eq:nablaXHXH}
    \nabla_{X_H} X_H = \pr^*\phi S - \mathsf{V}^\hor,
\end{equation}
\begin{equation}
    \label{eq:nablaXHS}
    \nabla_{X_H} S = \nabla_S X_H = \frac{1}{2} X_H,
\end{equation}
\begin{equation}
    \label{eq:nablaSS}
    \nabla_S S = \frac{1}{2} S,
\end{equation}
where we used the abbreviations
\begin{equation}
    \label{eq:tDef}
    t = \frac{1}{n+1} \Ric,
\end{equation}
\begin{equation}
    \label{eq:VDef}
    \mathsf{V} = \frac{2}{(n+1)(2n+1)} \mathsf{U},
\end{equation}
\begin{equation}
    \label{eq:phiDef}
    \phi = \frac{4}{(n+1)(2n+1)} \mathsf{f},
\end{equation}
and $\tau$ is the endomorphism corresponding to $t$ analogously to
\eqref{eq:varrhoDef}. In \cite{cahen.gutt.schwachhoefer:2005a}, the
following statement was obtained:
\begin{theorem}
    \label{theorem:FlatSymplectic}
    By \eqref{eq:nablaXhorYhor}--\eqref{eq:nablaSS} a flat symplectic
    torsion-free connection $\nabla$ is defined on $P$ and $\nabla$ is
    invariant under $X_H$ and $S$. Moreover, the third symmetrized
    covariant derivative of $H$ vanishes.
\end{theorem}
Recall that the operator of symmetrized covariant differentiation as a
derivation of the symmetric tensor product $\SymD:
\Gamma^\infty(\Sym^k T^*P) \longrightarrow \Gamma^\infty(\Sym^{k+1}
T^*P)$ is defined by
\begin{equation}
    \label{eq:SymDDef}
    (\SymD \gamma)(X_1, \ldots, X_{k+1})
    =
    \sum_{\ell=1}^{k+1} \left(\nabla_{X_\ell}\gamma\right)
    (X_1, \ldots, \stackrel{\ell}{\wedge}, \ldots, X_{k+1}),
\end{equation}
where $X_1, \ldots, X_{k+1} \in \Gamma^\infty(TP)$ and $\gamma \in
\Gamma^\infty(\Sym^k T^*P)$. Locally, $\SymD$ can be written as
\begin{equation}
    \label{eq:SymDLocally}
    \SymD = \D x^i \vee \nabla_{\frac{\partial}{\partial x^i}},
\end{equation}
where $\vee$ denotes the symmetrized tensor product.  Then
Theorem~\ref{theorem:FlatSymplectic} means
\begin{equation}
    \label{eq:DthreeHZero}
    \SymD^3 H = 0.
\end{equation}
\begin{remark}
    \label{remark:SymmetricRicciType}
    The Ricci type connection on $M$ is \emph{symmetric} iff
    $\mathsf{U} = 0$ in which case $\mathsf{f}$ turns out to be
    constant). This particular case has been studied in detail in
    \cite{bieliavsky.et.al:2006a}.
\end{remark}

%
%

\section{General Remarks on the Ricci Operator}
\label{sec:RicciOperator}

Before discussing the star products on $P$ and $M$, respectively, we
introduce the following second order differential operator on a
symplectic manifold with symplectic connection, which is also of
independent interest. Let $(M, \omega)$ be symplectic with a
torsion-free symplectic connection $\nabla$ (not necessarily of Ricci
type). Then the Ricci tensor $\Ric \in \Gamma^\infty(\Sym^2 T^*M)$ can
be used to define a `Laplace'-like operator $\DeltaRic$ as follows: We
denote by $\Ric^\sharp \in \Gamma^\infty(\Sym^2 TM)$ the symmetric
bivector obtained from $\Ric$ under the musical isomorphism $\sharp$
with respect to $\omega$.
\begin{definition}[Ricci operator]
    \label{definition:RicciOperator}
    The Ricci operator $\DeltaRic: C^\infty(M) \longrightarrow
    C^\infty(M)$ is defined by
    \begin{equation}
        \label{eq:DeltaRicDef}
        \DeltaRic u = \frac{1}{2} \SP{\Ric^\sharp, \SymD^2 u},
    \end{equation}
    where $\SP{\cdot, \cdot}$ denotes the natural pairing.
\end{definition}
If locally we write $\Ric^\sharp = \frac{1}{2} \Ric^{ij}
\frac{\partial}{\partial x^i} \vee \frac{\partial}{\partial x^j}$ then
\begin{equation}
    \label{eq:RicciOperatorLocally}
    \DeltaRic u = \Ric^{ij} 
    \left(\frac{\partial^2 u}{\partial x^i \partial x^j}
        - \Gamma^k_{ij} \frac{\partial u}{\partial x^k}
    \right),
\end{equation}
where $\Gamma^k_{ij}$ are the local Christoffel symbols of $\nabla$.

Since on a symplectic manifold we have a canonical volume form, the
Liouville form $\Omega$, one can ask whether $\DeltaRic$ is a
symmetric operator with respect to the $L^2$-inner product on
$C^\infty_0(M)$ induced by $\Omega$: in general, this is not the case.
However, there is an easy way to correct this. We need to recall some
basic features of the global symbol calculus for differential
operators on manifolds with connection, see
e.g.~\cite{bordemann.neumaier.pflaum.waldmann:2003a,
  bordemann.neumaier.waldmann:1998a}.

We denote by $(q, p)$ the canonical coordinates on $T^*U \subseteq
T^*M$ induced by a local chart $(U, x)$ of $M$. Then the
standard-ordered quantization of a function $f \in \Pol^\bullet(T^*M)$
on the cotangent bundle $\pi: T^*M \longrightarrow M$, which is
polynomial in the momenta, is defined by
\begin{equation}
    \label{eq:stdrepDef}
    \stdrep(f)u =
    \sum_{r=0}^\infty \frac{1}{r!} \left(\frac{\hbar}{\I}\right)^r 
    \frac{\partial^r f}{\partial p_{i_1} \cdots \partial p_{i_r}}
    \Big|_{p=0}
    \inss\left(\frac{\partial}{\partial x^{i_1}}\right)
    \cdots
    \inss\left(\frac{\partial}{\partial x^{i_r}}\right)
    \frac{1}{r!} \SymD^r u,
\end{equation}
where $u \in C^\infty(M)$ and $\inss$ denotes the (symmetric)
insertion map into the first argument. Clearly, \eqref{eq:stdrepDef}
is globally well-defined and does not depend on the coordinates. The
constant $\frac{\hbar}{\I}$ can safely be set to $1$ in our context;
however, we have included it for the sake of physical interpretation
of $\stdrep$ as a quantization map. Then \eqref{eq:stdrepDef} gives a
linear bijection $\stdrep: \Pol^\bullet(T^*M) \longrightarrow
\Diffop(M)$.

The symmetric algebra $\mathcal{S}^\bullet (TM) =
\bigoplus_{r=0}^\infty \Gamma^\infty(\Sym^r TM)$ is canonically
isomorphic to the polynomial functions $\Pol^\bullet (T^*M)$ as graded
associative algebra via the ``universal momentum map'' $\mathcal{J}$,
determined by $\mathcal{J}(u) = \pi^*u$ and
$(\mathcal{J}(X))(\alpha_q) = \alpha_q(X_(q))$ for $u \in
\mathcal{S}^0(TM) = C^\infty(M)$ and $X \in \mathcal{S}^1(TM) =
\Gamma^\infty(TM)$. Thus we can rephrase \eqref{eq:DeltaRicDef} as
\begin{equation}
    \label{eq:DeltaRicAsQuantization}
    \DeltaRic = - \frac{2}{\hbar^2} \stdrep(\mathcal{J}(\Ric^\sharp)),
\end{equation}
whence $\DeltaRic$ is the standard-ordered quantization of the
quadratic function $\mathcal{J}(\Ric^\sharp)$ on $T^*M$.

The standard-ordered quantization can be seen as a particular case of
the $\kappa$-ordered quantization which is obtained as follows. On
$C^\infty(T^*M)$ one has a Laplace operator $\Delta$ arising from the
pseudo-Riemannian metric $g_0$, which is defined by the natural
pairing of the vertical and horizontal (with respect to $\nabla$)
subspaces of $T(T^*M)$. Locally, $\Delta$ is given by
\begin{equation}
    \label{eq:DeltaIndefinite}
    \Delta\Big|_{T^*U}
    =
    \frac{\partial^2}{\partial q^i \partial p_i}
    +
    p_k \pi^*(\Gamma^k_{ij})
    \frac{\partial^2}{\partial p_i \partial p_j}
    +
    \pi^*(\Gamma^i_{ij}) \frac{\partial}{\partial p_j}.
\end{equation}
On polynomial functions $\Delta$ is just the covariant divergence
operator \cite[Eq.~(111)]{bordemann.neumaier.waldmann:1998a}, i.e.
\begin{equation}
    \label{eq:CovariantDivergence}
    \Delta \mathcal{J}(T) = \mathcal{J}(\divergence_\nabla T)
\end{equation}
for $T \in \mathcal{S}^k(TM)$ where with $\alpha_1, \ldots,
\alpha_{k-1} \in \Gamma^\infty(T^*M)$
\begin{equation}
    \label{eq:CovariantDivergenceDef}
    (\divergence_\nabla T) (\alpha_1, \ldots, \alpha_{k-1}) =
    \tr\left(
        X \mapsto (\nabla_X T)(\cdot, \alpha_1, \ldots, \alpha_{k-1})
    \right).
\end{equation}
Locally, $\divergence_\nabla = \inss(\D x^i)
\nabla_{\frac{\partial}{\partial x^i}}$.  Using $\Delta$, the
$\kappa$-ordered quantization is defined by
\cite{bordemann.neumaier.pflaum.waldmann:2003a}
\begin{equation}
    \label{eq:varrhokappa}
    \krep(f) = \stdrep\left(\E^{-\I\kappa\hbar\Delta} f\right),
\end{equation}
where in particular the \emph{Weyl-ordered case} $\weylrep =
\varrho_{\kappa = 1/2}$ is of interest for us. In general, we have for
the formal adjoint of $\krep(f)$ with respect to $\Omega$
\begin{equation}
    \label{eq:krepAdjoint}
    \krep(f)^\dagger 
    = \krep\left(\E^{-\I\hbar(1-2\kappa)\Delta} \cc{f}\right)
\end{equation}
whence $\weylrep(f)^\dagger = \weylrep(\cc{f})$. Thus
$\weylrep(\mathcal{J}(\Ric^\sharp))$ gives a symmetric
operator. Explicitly, one finds using \eqref{eq:CovariantDivergence}
\begin{equation}
    \label{eq:kappaOrderedRiccian}
    -\frac{2}{\hbar^2} \krep(\mathcal{J}(\Ric^\sharp)) =
    \DeltaRic + 2\kappa \Lie_{\divergence_\nabla \Ric^\sharp} +
    \kappa^2 \divergence_\nabla^2(\Ric^\sharp),
\end{equation}
since no higher order terms contribute thanks to $\Ric^\sharp \in
\mathcal{S}^2(TM)$. Moreover, $\kappa^2 \divergence_\nabla^2
(\Ric^\sharp)$ is already a multiplication operator with a \emph{real}
function and hence symmetric itself for all $\kappa$. For $\kappa =
\frac{1}{2}$ we have:
\begin{lemma}
    \label{lemma:SymmetricRiccian}
    The operator $\DeltaRic + \Lie_{\divergence_\nabla \Ric^\sharp} +
    \frac{1}{4} \divergence_\nabla^2 \Ric^\sharp$ as well as the operator
    $\DeltaRic + \Lie_{\divergence_\nabla \Ric^\sharp}$ are symmetric.
\end{lemma}

We conclude this section with the computation of the covariant
divergences of $\Ric^\sharp$ in the case of a Ricci type
connection. In view of equation  \eqref{eq:nablaXvarrho} and the
fact that $\nabla_X$ commutes with $\sharp$, we have :
\begin{equation}
    \label{eq:nablaXRicsharp}
    \nabla_X \Ric^\sharp 
    = \left(\nabla_X \Ric\right)^\sharp
    = \frac{1}{2n+1}\left(X^\flat \vee \mathsf{U}^\flat\right)^\sharp
    = \frac{1}{2n+1} X \vee \mathsf{U}
\end{equation}
whence
\begin{equation}
    \label{eq:DivergenceRic}
    \divergence_\nabla \Ric^\sharp = \mathsf{U}.
\end{equation}
Moreover, from \eqref{eq:nablaXU} and \eqref{eq:tracevarrho} we get
\begin{equation}
    \label{eq:DoubleDivRic}
    \divergence_\nabla^2 \left(\Ric^\sharp\right)
    = \divergence_\nabla \mathsf{U}
    = - \frac{2n +1}{2(n+1)} K + 2(n+1) f.
\end{equation}
We can also obtain from this the symmetric version of the Ricci operator,
 according to Lemma~\ref{lemma:SymmetricRiccian}.

%
%

\section{The Weyl-Moyal Star Product}
\label{sec:WeylMoyalStarProduct}

After this excursion on the Ricci operator, we are now back to the
situation of Section~\ref{sec:RicciTypeConnections}. On the big
space $P$, the
symplectic connection $\nabla$ is flat. We have thus a Weyl-Moyal star
product on $P$, explicitly given by
\begin{equation}
    \label{eq:WeylMoyal}
    f \star g = \sum_{r=0}^\infty \frac{1}{r!}
    \left(\frac{\nu}{2}\right)^r C_r(f, g)
\end{equation}
for $f, g \in C^\infty(P)[[\nu]]$, where
\begin{equation}
    \label{eq:CrExplicit}
    C_r(f, g) = 
    \SP{\Lambda_P \otimes \cdots \otimes \Lambda_P,
      \frac{1}{r!} \SymD^r f \otimes \frac{1}{r!} \SymD^r g}
\end{equation}
and the natural pairing is done ``over cross''. Locally we have
\begin{align}
    \nonumber
    C_r(f, g) &= 
    \Lambda_P^{i_1j_1} \cdots \Lambda_P^{i_rj_r}
    \inss\left(\frac{\partial}{\partial x^{i_1}}\right) \cdots
    \inss\left(\frac{\partial}{\partial x^{i_r}}\right) 
    \frac{1}{r!} \SymD^r f \\
    &\quad {\times} \;
    \inss\left(\frac{\partial}{\partial x^{j_1}}\right) \cdots
    \inss\left(\frac{\partial}{\partial x^{j_r}}\right) 
    \frac{1}{r!} \SymD^r g.
    \label{eq:CrLocally}
\end{align}
Since $\nabla$ is flat  $\star$ defines  an associative law on the space of formal power series in the parameter $\nu$ with coefficients in $C^\infty(P)$  \cite{bayen.et.al:1978a}. The $C_r$ are bidifferential operators which are of order at most $r$ in each argument, and satisfy the symmetry condition
$C_r(u,v)=(-1)^rC_r(v,u)$.  These properties are summarized by saying that $\star$
is a  natural star product of Weyl-type as.

We have now two derivations of $\star$: first it follows directly from
the fact that $\nabla$ is $S$-invariant and $\Lie_S \Lambda_P = -
\Lambda_P$ that
\begin{equation}
    \label{eq:nuEuler}
    \mathcal{E} = \nu \frac{\partial}{\partial \nu} + \Lie_S
\end{equation}
is a \emph{$\nu$-Euler derivation}, i.e. a $\mathbb{C}$-linear
derivation of $\star$, see e.g.\cite{gutt.rawnsley:1999a}. Second, we
consider the quasi-inner derivation $\frac{1}{\nu} \ad(H)$. Since
$\SymD^3 H = 0$, only the terms of order $\nu^1$ and $\nu^2$
contribute. But thanks to the Weyl-type symmetry of $\star$, only odd
powers of $\nu$ occur in commutators, this is immediate from
\eqref{eq:CrLocally}. Thus we have
\begin{equation}
    \label{eq:adH}
    \frac{1}{\nu} \ad(H)f = \frac{1}{\nu} \nu\{H, f\} = X_H(f),
\end{equation}
which shows that $X_H$ is a quasi-inner derivation. With other words,
$\star$ is \emph{strongly invariant} with respect to the (classical)
momentum map $H$. This strong invariance will allow us to use the
phase space reduction also for $\star$ to obtain a star product on
$M$. To this end, we first note the following:
\begin{lemma}
    \label{lemma:pruStarprv}
    There are unique bidifferential operators $\hat{C}^\red_r$ on $M$ of
    order $r$ in each argument such that
    \begin{equation}
        \label{eq:CredDef}
        \pr^*u \star \pr^*v 
        = \sum_{r=0}^\infty \frac{1}{r!} \left(\frac{\nu}{2}\right)^r
        \frac{1}{H^r} \pr^* \hat{C}_r^\red(u, v).
    \end{equation}
    In particular, we have $\hat{C}_1^\red(u, v) = \SP{\Lambda_M, \D u
      \otimes \D v}$.
\end{lemma}
\begin{proof}
    Clearly, $H^r C_r (\pr^*u, \pr^*v)$ is invariant under $X_H$ and
    homogeneous of degree $0$ under $S$, hence a pull-back of a
    function $\hat{C}_r^\red(u, v) \in C^\infty(M)$ via $\pr$. This
    defines $\hat{C}_r^\red$ uniquely. The statement on the order of
    differentiation is straightforward. Finally, $\hat{C}_1^\red$ is
    obtained from \eqref{eq:PoissonBracketuv}.
    \qed
\end{proof}
\begin{remark}
    \label{remark:TooNaive}
    Though it seems tempting, the operators $\hat{C}_r^\red$ do
    \emph{not} combine into a star product on $M$ directly: the
    prefactors $H^r$ spoil the associativity as one can show by a
    direct computation. Hence we will need a slightly more involved
    reduction.
\end{remark}

We will need the second order term of $H \star f$ with an arbitrary
function $f \in C^\infty(P)$. By the symmetry properties of
\eqref{eq:CrLocally} we know that $C_2(H, f) = C_2(f, H)$.
\begin{proposition}
    \label{proposition:CtwoExplicitly}
    Let $f \in C^\infty(P)$, then
    \begin{equation}
        \label{eq:CtwoExplicit}
        \begin{split}
            C_2(f, H) 
            = 
            \frac{1}{H}\Bigg(
            &- \frac{1}{n+1} \left(\DeltaRic\right)^\hor f 
            - \frac{1}{2} \Lie_{\mathsf{V}^\hor} f
            - \pr^*(\phi) \Lie_S^{\:2} f
            \\
            &+ \left(
                \Lie_{\mathsf{V}^\hor} \Lie_S 
                + \Lie_S \Lie_{\mathsf{V}^\hor}
            \right) f
            + \frac{1}{2} \Lie_{X_H}^{\:2} f
            \Bigg).
        \end{split}
    \end{equation}
    where $ \left(\DeltaRic\right)^\hor f $ is the pairing of the
    horizontal lift of the $2$-tensor $\Ric^{M\sharp}$ with the second
    covariant derivative with respect to the flat connection $\nabla$
    on $P$. Precisely, in a local chart, it is given by
    \begin{equation}
        \label{eq:RicciOperatorHor}
        \pr^*(\Ric^{M\sharp ij})
        [\Lie_{\partial_{i}^\hor}\Lie_{\partial_{j}^\hor} f
        -  \pr^*(\Gamma^{Mk}_{ij}) \Lie_{\partial_{k}^\hor} f
        + \frac{1}{n+1} \pr^*  \Ric^M_{ij} \Lie_S f].
    \end{equation}
    Observe that $\Ric^{M\sharp ij}\Ric^M_{ij}=-\tr(\rho^2)$.
\end{proposition}
\begin{proof}
    First we note that $\inss(X)\SymD^2 f = 2 \nabla_X \D f$ by the
    very definition of $\SymD$. Thus using \eqref{eq:PoissonP} we can
    compute $C_2(f, g)$ for $f, g \in C^\infty(P)$ explicitly and get
    \begin{align}
        C_2(f, g)
        &= \frac{1}{H^2}
        \Bigg(
        \pr^*\left(\Lambda_\red^{ij} \Lambda_\red^{kl}\right)
        \left(\nabla_{\partial_k^\hor}\D f\right)(\partial_i^\hor)
        \left(\nabla_{\partial_l^\hor}\D g\right)(\partial_j^\hor)
        \nonumber\\
        &+
        \pr^*\left(\Lambda_\red^{kl}\right)
        \left(\nabla_{\partial_k^\hor}\D f\right)(S)
        \left(\nabla_{\partial_l^\hor}\D g\right)(X_H) 
        \nonumber\\
        &+
        \pr^*\left(\Lambda_\red^{kl}\right)
        \left(\nabla_{\partial_l^\hor}\D f\right)(X_H)
        \left(\nabla_{\partial_k^\hor}\D g\right)(S) 
        \nonumber\\
        &+
        \pr^*\left(\Lambda_\red^{ij}\right)
        \left(\nabla_S \D f\right)(\partial_i^\hor)
        \left(\nabla_{X_H} \D g\right)(\partial_j^\hor) 
        \nonumber\\
        &+
        \pr^*\left(\Lambda_\red^{ij}\right)
        \left(\nabla_{X_H} \D f\right)(\partial_j^\hor)
        \left(\nabla_S \D g\right)(\partial_i^\hor)
        \nonumber\\
        &+
        \left(\nabla_S \D f\right)(S)
        \left(\nabla_{X_H} \D g\right)(X_H)
        +
        \left(\nabla_{X_H} \D f\right)(X_H)
        \left(\nabla_S \D g\right)(S)
        \nonumber\\
        &-
        \left(\nabla_S \D f\right)(X_H)
        \left(\nabla_{X_H} \D g\right)(S) 
        -
        \left(\nabla_{X_H} \D f\right)(S)
        \left(\nabla_S \D g\right)(X_H)
        \Bigg),
        \label{eq:Czwei}
    \end{align}
    where we have used local coordinates on $M$ as well as the
    horizontal lift according to \eqref{eq:horizontalLift}. Next, one
    computes the second covariant derivatives of $H$ explicitly. One
    finds
    \begin{equation}
        \label{eq:nablaXhordHYhor}
        \left(\nabla_{X^\hor} \D H\right) \left(Y^\hor\right) 
        = - \pr^*(t(X, Y)) H,
    \end{equation}
    \begin{equation}
        \label{eq:nablaXhordHXH}
        \left(\nabla_{X^\hor} \D H\right) \left(X_H\right) 
        =
        \pr^*(\omega(X, \mathsf{V})) H
        =
        \left(\nabla_{X_H} \D H\right) \left(X^\hor\right),
    \end{equation}
    \begin{equation}
        \label{eq:nablaXhordHS}
        \left(\nabla_{X^\hor} \D H\right) \left(S\right) 
        =
        0
        =
        \left(\nabla_S \D H\right) \left(X^\hor\right),
    \end{equation}
    \begin{equation}
        \label{eq:nablaXHdHS}
        \left(\nabla_{X_H} \D H\right) \left(S\right) 
        =
        0
        =
        \left(\nabla_{S} \D H\right) \left(X_H\right),
    \end{equation}
    \begin{equation}
        \label{eq:nablaXHdHXHnablasdHS}
        \left(\nabla_{X_H} \D H\right) \left(X_H\right) 
        =
        - \pr^*(\phi) H
        \quad
        \textrm{and}
        \quad
        \left(\nabla_{S} \D H\right) \left(S\right)
        =
        \frac{1}{2} H,
    \end{equation}
    where we used $\D H(X_H) = 0$ and $\D H(S) = H$ as well as $\D
    H(X^\hor) = 0$ together with the explicit formulas
    \eqref{eq:nablaXhorYhor} -- \eqref{eq:nablaSS}. Putting things
    together gives the result thanks to the local form
    \eqref{eq:RicciOperatorLocally} for $\DeltaRic$ and $\Ric^{ij} =
    \Lambda_M^{ik} \Lambda_M^{jl} \Ric_{kl}$.
    \qed
\end{proof}
\begin{definition}
    \label{definition:Delta}
    For $f \in C^\infty(P)$ we define the second order differential
    operator $\Delta$ by
    \begin{equation}
        \label{eq:DeltaDef}
        \Delta f = C_2(f, H).
    \end{equation}
\end{definition}
As we shall see in the next section, this operator will be crucial for
the construction of the reduced star product.

%
%

\section{Reduction of the Star Product}
\label{sec:ReductionStarProducts}

We come now to the reduction of $\star$. Here we follow essentially
the BRST / Koszul approach advocated in
\cite{bordemann.herbig.waldmann:2000a} which simplifies drastically
thanks to the codimension one reduction. Codimension one reductions
have also been discussed by Glößner \cite{gloessner:1998a:pre,
  gloessner:1998a} and Fedosov\cite{fedosov:1994b}, while the general
case of reduction for coisotropic constraint manifolds is discussed in
\cite{bordemann:2005a, bordemann:2004a:pre, cattaneo.felder:2005a:pre,
  cattaneo.felder:2004a}. We shall briefly recall those aspects of
\cite{bordemann.herbig.waldmann:2000a} which are needed here.

We first consider the classical situation: the classical Koszul
operator $\partial: C^\infty(P) \longrightarrow C^\infty(P)$ is
defined by
\begin{equation}
    \label{eq:ClassicalKoszul}
    \partial f = f (H-1).
\end{equation}
Next we define the classical homotopy $h: C^\infty(P) \longrightarrow
C^\infty(P)$ by
\begin{equation}
    \label{eq:ClassicalHomotopy}
    h f =
    \begin{cases}
        \frac{1}{H-1} \left(f - \pr_1^* \iota^* f\right)
        & \textrm{on} \; P \setminus \iota(\Sigma) \\
        \Lie_S f
        & \textrm{on} \; \iota(\Sigma).
    \end{cases}
\end{equation}
An easy argument shows that $hf$ is actually smooth. Then we have the
homotopy formula
\begin{equation}
    \label{eq:ClassicalExactness}
    f = \partial h f + \pr_1^* \iota^* f
\end{equation}
for $f \in C^\infty(P)$ together with the properties
\begin{equation}
    \label{eq:iotapartialNullhprNull}
    \iota^* \partial = 0
    \quad
    \textrm{and}
    \quad
    h \pr_1^* = 0.
\end{equation}
Moreover, $\partial$, $h$, $\iota^*$, and $\pr_1^*$ are equivariant
with respect to the action of $X_H$ on $P$ and the Reeb vector field
on $\Sigma$, respectively. The classical \emph{vanishing ideal}
$\mathcal{I}_\Sigma$ of $\Sigma$ is given by
\begin{equation}
    \label{eq:ClassicalJSigma}
    \mathcal{I}_\Sigma = \ker \iota^* = \image \partial
\end{equation}
by \eqref{eq:ClassicalExactness}, and turns out to be a Poisson
subalgebra of $C^\infty(P)$ as $\Sigma$ is coisotropic. Moreover, 
\begin{equation}
    \label{eq:BSigmaDef}
    \mathcal{B}_\Sigma = 
    \left\{f
        \in C^\infty(P) 
        \; \big| \; 
        \{f, \mathcal{I}_\Sigma\} \subseteq \mathcal{I}_\Sigma 
    \right\}
\end{equation}
is the largest Poisson subalgebra of $C^\infty(P)$ such that
$\mathcal{I}_\Sigma \subseteq \mathcal{B}_\Sigma$ is a Poisson
ideal. It is well-known and easy to see that $\mathcal{B}_\Sigma$ can
also be characterized by
\begin{equation}
    \label{eq:BSigmaInvariance}
    \mathcal{B}_\Sigma =
    \left\{
        f \in C^\infty(P)
        \; \big| \; 
        \iota^*f \in C^\infty(\Sigma)^Z = \pi^* C^\infty(M)
    \right\}
\end{equation}
from which it easily follows that the Poisson algebra
$\mathcal{B}_\Sigma \big/ \mathcal{I}_\Sigma$ is isomorphic to the
Poisson algebra $C^\infty(M)$ via $\iota^*$ and $\pi^*$.

We shall now deform the above picture according to
\cite{bordemann.herbig.waldmann:2000a} where we only use the `Koszul
part' of the BRST complex. First we define the quantum Koszul operator
$\qkoszul: C^\infty(P)[[\nu]] \longrightarrow C^\infty(P)[[\nu]]$ by
\begin{equation}
    \label{eq:QuantumKoszul}
    \qkoszul f = f \star (H-1).
\end{equation}
We set
\begin{equation}
    \label{eq:ImageqKoszul}
    \qI = \image \qkoszul = 
    \left\{
        f \in C^\infty(P)[[\nu]]
        \; \big| \; 
        f = g \star (H-1)
    \right\},
\end{equation}
which is the left ideal generated by $H-1$ with respect to
$\star$. Next we consider
\begin{equation}
    \label{eq:qBDef}
    \qB =
    \left\{
        f \in C^\infty(P)[[\nu]]
        \; \big| \;
        [f, \qI]_\star \subseteq \qI
    \right\},
\end{equation}
which is the largest subalgebra of $C^\infty(P)[[\nu]]$ such that $\qI
\subseteq \qB$ is a two-sided ideal, the so-called \emph{idealizer} of
$\qI$. The following simple algebraic lemma is at the core of
Bordemann's interpretation of the reduction procedure
\cite{bordemann:2004a:pre, bordemann:2005a}:
\begin{lemma}
    \label{lemma:AredDef}
    Let $\mathcal{A}$ be a unital $\mathbb{k}$-algebra and
    $\mathcal{I} \subseteq \mathcal{A}$ a left ideal. Let $\mathcal{B}
    \subseteq \mathcal{A}$ be the idealizer of $\mathcal{I}$ and
    $\mathcal{A}_\red = \mathcal{B} \big/\mathcal{I}$. Then the left
    $\mathcal{A}$-module $\mathcal{A} \big/\mathcal{I}$ becomes a
    right $\mathcal{A}_\red$-module via $[b]: [a] \mapsto [ab]$ for
    $[b] \in \mathcal{A}_\red$ and $[a] \in \mathcal{A} \big/
    \mathcal{I}$, such that
    \begin{equation}
        \label{eq:AredIsomorphicEndoAI}
        \mathcal{A}_\red 
        \cong 
        \End_{\mathcal{A}} 
        \left(\mathcal{A}\big/\mathcal{I}\right)^\opp.
    \end{equation}
    This way $\mathcal{A} \big/\mathcal{I}$ becomes a $(\mathcal{A},
    \mathcal{A}_\red)$-bimodule.
\end{lemma}

In our situation, we want to show that $C^\infty(P)[[\nu]] \big/ \qI$
provides a deformation of $C^\infty(\Sigma)$ and $\qB \big/ \qI$
induces a star product $\starred$ on $M$. This can be done in great
generality, in our situation the arguments simplify thanks to the
codimension one case.

We define the quantum homotopy
\begin{equation}
    \label{eq:qhDef}
    \qh =  h \left(\id + (\qkoszul - \partial) h\right)^{-1}
\end{equation}
and the quantum restriction map
\begin{equation}
    \label{eq:qiotaDef}
    \qiota = \iota^*\left(\id + (\qkoszul - \partial) h\right)^{-1},
\end{equation}
which are clearly well-defined since $\qkoszul - \partial$ is at least
of order $\lambda$ and thus $\id + (\qkoszul - \partial) h$ is
invertible by a geometric series. From \eqref{eq:ClassicalExactness}
we immediately find
\begin{equation}
    \label{eq:QuantumExactness}
    f = \qkoszul \qh f + \pr_1^* \qiota f
\end{equation}
for all $f \in C^\infty(P)[[\nu]]$. Moreover, we still have the
relations
\begin{equation}
    \label{eq:qiotapridqhprNull}
    \qiota \pr_1^* = \id_{C^\infty(\Sigma)[[\nu]]}
    \quad
    \textrm{and}
    \quad
    \qh \pr_1^* = 0,
\end{equation}
as in the classical case. Finally, all maps are still equivariant with
respect to the flow of $X_H$ and $Z$, respectively.
\begin{proposition}
    \label{proposition:LeftModule}
    The functions $C^\infty(\Sigma)[[\nu]]$ becomes a left
    $\star$-module via
    \begin{equation}
        \label{eq:LeftModule}
        f \bullet \psi = \qiota(f \star \pr_1^* \psi).
    \end{equation}
    This module structure is isomorphic to the module structure of
    $C^\infty(P)[[\nu]] \big/ \qI$ via $\pr_1^*$ and $\qiota$.
\end{proposition}
\begin{proof}
    Since $\image \qkoszul = \ker \qiota$ by
    \eqref{eq:QuantumExactness}, the quantum restriction map $\qiota$
    induces a linear bijection $C^\infty(P)[[\nu]] \big/ \qI
    \longrightarrow C^\infty(\Sigma)[[\nu]]$ whose inverse is induced
    by $\pr_1^*$ by \eqref{eq:qiotapridqhprNull}. Then
    \eqref{eq:LeftModule} is just the pulled-back module structure.
    \qed
\end{proof}

From the strong invariance of the star product $\star$ we obtain the
following characterization of $\qB$:
\begin{lemma}
    \label{lemma:qBAlternative}
    $\qB = \left\{f \in C^\infty(P)[[\nu]] \; \big| \; \qiota f \in
        \pi^*C^\infty(M)[[\nu]] \right\}$.
\end{lemma}
\begin{proof}
    Indeed, on one hand we have $Z \qiota f = - \nu \qiota \ad(H) f$
    by equivariance of $\qiota$. Since for $f \in \qB$ we have $\ad(H)
    f \in \qI$, this gives $\qiota f = 0$ by
    \eqref{eq:QuantumExactness}. Conversely, $\qiota f = 0$ implies
    $\ad(H) f \in \qI$ and hence $(H - 1) \star f \in \qI$ from which
    $\qI \star f \subseteq \qI$ follows. But this implies $f \in \qB$.
    \qed
\end{proof}
\begin{theorem}
    \label{theorem:QuotientAlgebra}
    The quotient $\qB \big/ \qI$ is $\mathbb{C}[[\nu]]$-linearly
    isomorphic to $C^\infty(M)[[\nu]]$ via
    \begin{equation}
        \label{eq:MtoBmodI}
        C^\infty(M)[[\nu]] \ni u 
        \; \mapsto \; 
        [\pr^*u] \in \qB \big/ \qI
    \end{equation}
    with inverse induced by
    \begin{equation}
        \label{eq:qiotaIso}
        \qB \big/ \qI \ni [f]
        \; \mapsto \;
        \qiota f \in \pi^*C^\infty(M)[[\nu]].
    \end{equation}
    This induces a deformed product $\starred$ for
    $C^\infty(M)[[\nu]]$ via
    \begin{equation}
        \label{eq:starredDef}
        \pi^*(u \starred v) = \qiota(\pr^*u \star \pr^*v),
    \end{equation}
    which turns out to be a differential star product quantizing
    $\{\cdot, \cdot\}_M$. Finally, the bimodule structure of
    $C^\infty(\Sigma)[[\nu]]$ according to Lemma~\ref{lemma:AredDef}
    and Proposition~\ref{proposition:LeftModule} is given by
    \begin{equation}
        \label{eq:RightModule}
        \psi \bullet u = \qiota(\pr_1^* \psi \star \pr^* u).
    \end{equation}
\end{theorem}
\begin{proof}
    For $u \in C^\infty(M)[[\nu]]$ we clearly have $\qiota \pr^* u =
    \qiota \pr_1^* \pi^* u = \pi^*u \in \pi^*C^\infty(M)[[\nu]]$
    whence $\pr^* u \in \qB$ by Lemma~\ref{lemma:qBAlternative} and
    \eqref{eq:MtoBmodI} is well-defined. Since $\qI = \ker \qiota$, by
    Lemma~\ref{lemma:qBAlternative} it follows that
    \eqref{eq:qiotaIso} is well-defined and injective. Clearly,
    \eqref{eq:MtoBmodI} and \eqref{eq:qiotaIso} are mutually inverse
    by \eqref{eq:qiotapridqhprNull} whence $\starred$ is an
    associative $\mathbb{C}[[\nu]]$-bilinear product for
    $C^\infty(M)[[\nu]]$. It can be shown
    \cite[Lem.~27]{bordemann.herbig.waldmann:2000a} that there exists
    a formal series $S = \id + \sum_{r=1}^\infty \nu^r S_r$ of
    differential operators $S_r$ on $P$ such that $\qiota = \iota^*
    \circ S$, from which it easily follows that $\starred$ is
    bidifferential. Finally, computing the first orders of
    \begin{equation}
        \label{eq:starredByOrders}
        u \starred v 
        = \sum_{r=0}^\infty \frac{1}{r!} 
        \left(\frac{\nu}{2}\right)^r
        C_r^\red(u, v)
    \end{equation}
    explicitly gives $C_0^\red(u, v) = uv$ and $C_1^\red(u, v) = \{u,
    v\}_M$, whence $\starred$ is a star product quantizing the correct
    Poisson bracket. The last statement is clear by construction.
    \qed
\end{proof}

Up to now we just followed the general reduction scheme from
\cite{bordemann.herbig.waldmann:2000a} which simplifies for the
codimension one case, see also \cite{gloessner:1998a:pre,
  gloessner:1998a}. Let us now bring our more specific features into
the game:
\begin{lemma}
    \label{lemma:DeltaAgain}
    Let $f \in C^\infty(P)^{X_H}[[\nu]]$ be invariant. Then
    \begin{equation}
        \label{eq:DeltaAgain}
        (\qkoszul - \partial)f = \frac{\nu^2}{8} \Delta f.
    \end{equation}
    The operator $\Delta$ is invariant.
\end{lemma}
\begin{proof}
    By invariance of $f$ we only have the second order term in the
    right multiplication by $H-1$, which was computed in
    Proposition~\ref{proposition:CtwoExplicitly}. The invariance of
    $\Delta$ follows from the strong invariance of $\star$ by
    \begin{align*}
        X_H(\Delta f) 
        &= X_H\left(f \star H - \frac{\nu}{2} C_1(f, H) - fH\right) \\
        &= \frac{1}{\nu} \ad(H)(f \star H) 
        + \frac{\nu}{2} X_H(X_H(f)) - X_H(f) H \\
        &= \frac{1}{\nu} (\ad(H)(f)) \star H
        - \frac{\nu}{2} C_1(X_H(f), H) - X_H(f) H \\
        &= \Delta(X_H(f)).
        \qed
    \end{align*}
\end{proof}
\begin{lemma}
    \label{lemma:qiotaExplicit}
    Let $f \in C^\infty(P)^{X_H}[[\nu]]$ be invariant. Then
    \begin{equation}
        \label{eq:qiotaExplicit}
        \qiota f 
        = \iota^*\left(\id + \frac{\nu^2}{8} \Delta h \right)^{-1} f.
    \end{equation}
\end{lemma}
\begin{proof}
    Since $h$ and $\Delta$ preserve invariance this follows by
    induction from the last lemma. 
    \qed
\end{proof}
\begin{proposition}
    \label{proposition:starredExplicit}
    Let $u, v \in C^\infty(M)[[\nu]]$. Then $u \starred v$ is
    determined by
    \begin{equation}
        \label{eq:ustarredv}
        \pi^*(u \starred v)
        =
        \iota^*\left(
            \sum_{r=0}^\infty \nu^r \sum_{2s + t = r}
            \frac{1}{(-8)^s 2^t t!}
            \left(\Delta h\right)^s
            \left(\frac{1}{H^t} \pr^*\hat{C}^\red_t(u, v)\right)
        \right).
    \end{equation}
\end{proposition}
\begin{proof}
    This is a simple consequence of Lemma~\ref{lemma:qiotaExplicit}
    together with Lemma~\ref{lemma:pruStarprv}.
    \qed
\end{proof}

From this formula we see that we have to proceed in two steps in order
to compute the true bidifferential operators $C_r^\red$ out of the
operators $\hat{C}^\red_r$: first we have to control $\Delta$ and $h$
applied to functions of the form $\frac{1}{H^t} \pr^*u$. Second, the
application of $\iota^*$ simply sets $H = 1$ and gives $\iota^* \pr^*
= \pi^*$, whence this part can be considered to be trivial.
\begin{lemma}
    \label{lemma:HomotopyHtpr}
    Let $u \in C^\infty(M)$ and $k \in \mathbb{N}$. Then
    \begin{equation}
        \label{eq:hHtpr}
        h \left(\frac{1}{H^k} \pr^*u\right)
        =
        - \left(\frac{1}{H} + \cdots + \frac{1}{H^k}\right) \pr^* u.
    \end{equation}
\end{lemma}
\begin{proof}
    On $P \setminus \iota(\Sigma)$ we have from
    \eqref{eq:ClassicalHomotopy}
    \begin{align*}
        h \left(\frac{1}{H^k} \pr^*u\right)
        &=
        \frac{1}{H-1} \left(
            \frac{1}{H^k} \pr^*u 
            - \pr_1^*\iota^*\left(\frac{1}{H^k} \pr^*u\right)
        \right) \\
        &=
        \frac{1}{H-1} \left(\frac{1}{H^k} - 1\right) \pr^*u \\
        &= - \left(\frac{1}{H} + \cdots + \frac{1}{H^k}\right) \pr^* u.
    \end{align*}
    By continuity this extends also to $\iota(\Sigma)$.
    \qed
\end{proof}

In particular, applying the homotopy $h$ to a linear combination of
functions of the form $\frac{1}{H^k}\pr^*u$ gives again such a linear
combination.
\begin{lemma}
    \label{lemma:DeltaHkpru}
    Let $u \in C^\infty(M)$ and $k \in \mathbb{N}$. Then
    \begin{equation}
        \label{eq:DeltaHkpru}
        \begin{split}
            \Delta \left(
                \frac{1}{H^k} \pr^*u
            \right)
            &= 
            - \frac{1}{n+1} \frac{1}{H^{k+1}} 
            \pr^*\bigg(
            \DeltaRic u \\
            &\qquad
            + \frac{4k+1}{2n+1} \Lie_{\mathsf{U}} u
            + \frac{k}{n+1} \tr(\varrho^2) u
            + \frac{4k^2}{2n+1} \mathsf{f} u
            \bigg).
        \end{split}
    \end{equation}
\end{lemma}
\begin{proof}
    This follows from $\Lie_{X_H} H = 0 = \Lie_{X^\hor} H$ and $\Lie_S
    H = H$ together with $\Lie_{X^\hor} \pr^* u = \pr^* \Lie_X u$ and
    $\Lie_{X_H} \pr^*u = 0 = \Lie_S \pr^*u$ as well as
    Proposition~\ref{proposition:CtwoExplicitly}.
    \qed
\end{proof}

Again, we see that applying $\Delta$ to this particular class of
functions reproduces such linear combinations, though, of course, the
combinatorics gets involved. In contrast to \eqref{eq:hHtpr}, the
function $u$ is differentiated in \eqref{eq:DeltaHkpru}.
\begin{remark}
    \label{remark:ReallyExplicit}
    From the above two lemmas the star product $\starred$ can be
    computed by Proposition~\ref{proposition:starredExplicit} in all
    orders. However, the explicit evaluation of the iteration $(\Delta
    h)^s$ seems to be tricky: the combinatorics gets quite involved,
    even in the case, where $\nabla^M$ is symmetric, i.e $\mathsf{U} =
    0$ and $\mathsf{f} = \const$, $\tr(\varrho^2) = \const$. In this
    particular case only the Ricci operator $\DeltaRic$ has to be
    applied successively to the $\hat{C}^\red_r$ to produce the
    $C^\red_r$, including, of course, still some non-trivial
    combinatorics.
\end{remark}

%
%

\section{Properties of the Reduced Star Product}
\label{sec:Properties}

In this section we collect some further properties of the reduced star
product.
\begin{proposition}
    \label{proposition:WeylType}
    The star product $\starred$ is of Weyl-type, i.e. the
    bidifferential operators $C^\red_r$ satisfy
    \begin{equation}
        \label{eq:ccCredr}
        \cc{C^\red_r(u, v)} = C^\red_r(\cc{u}, \cc{v})
    \end{equation}
    and
    \begin{equation}
        \label{eq:AntisymmCredr}
        C^\red_r(u, v) = (-1)^r C^\red_r(v, u).
    \end{equation}
    In particular, complex conjugation becomes a $^*$-involution
    \begin{equation}
        \label{eq:ccStarInv}
        \cc{u \starred v} = \cc{v} \starred \cc{u},
    \end{equation}
    where $\cc{\nu} = - \nu$ by definition. Moreover, $\starred$ is
    natural in the sense of \cite{gutt.rawnsley:2003a}, i.e.
    $C^\red_r$ is a bidifferential operator of order $r$ in each
    argument.
\end{proposition}
\begin{proof}
    The operators $\hat{C}^\red_r$ clearly satisfy \eqref{eq:ccCredr}
    and \eqref{eq:AntisymmCredr}. From \eqref{eq:starredByOrders} it
    follows that the operators $C^\red_r$ are real operators, i.e.
    they satisfy~\eqref{eq:ccCredr}, since $\Delta$ and $h$ are real.
    Moreover, by \eqref{eq:ustarredv}, only the operators
    $\hat{C}^\red_{r-2s}$ with $0 \le 2s \le r$ contribute to
    $C^\red_r$ whence $C^\red_r$ also satisfy \eqref{eq:AntisymmCredr}
    as $2s$ is even. Then \eqref{eq:ccStarInv} follows. Moreover,
    $\hat{C}^\red_r$ is a differential operator of order $r$ in each
    argument. Since the application of $\Delta h$ to $\frac{1}{H^k}
    \pr^*u$ by Lemma~\ref{lemma:HomotopyHtpr} and
    Lemma~\ref{lemma:DeltaHkpru} gives a second order differential
    operator on the part $\pr^*u$, we conclude that $C^\red_r$ is of
    order $r$ in each argument, too.
    \qed
\end{proof}

In a next step, we explicitly compute the second order term $C^\red_2$
of $\starred$. According to \cite{gutt.rawnsley:2003a}, the second
order term of a natural star product determines uniquely a symplectic
connection: in our case, we reproduce the Ricci type connection
$\nablaM$.
\begin{proposition}
    \label{proposition:SecondOrderTerm}
    Let $u, v \in C^\infty(M)[[\nu]]$. Then
    \begin{equation}
        \label{eq:CTwoRed}
        C^\red_2(u, v)
        =
        \SP{\Lambda_\red \otimes \Lambda_\red,
          \frac{1}{2}\SymD_M^2 u \otimes \frac{1}{2}\SymD_M^2 v}
        +
          \frac{2}{n+1}\SP{\Ric^\sharp, \D u \otimes \D v}.
    \end{equation}
    In particular, the symplectic connection determined by $\starred$
    is $\nablaM$.
\end{proposition}
\begin{proof}
    From \eqref{eq:ustarredv} we see that
    \begin{align*}
        \pi^* C^\red_2(u, v)
        &= 8\iota^*\left(
            \sum_{2s+t=2} \frac{1}{(-8)^s 2^t t!} (\Delta h)^s
            \left(\frac{1}{H^t} \pr^* \hat{C}^\red_2 (u, v)\right) 
        \right)\\
        &= \iota^*\left(
            \frac{1}{H^2} \pr^* \hat{C}^\red_2 (u, v)
            - \Delta h \pr^* \hat{C}^\red_0(u, v)
        \right) \\
        &= \pi^* \hat{C}^\red_2 (u, v),
    \end{align*}
    since $h \pr^* = 0$ by the very definition
    \eqref{eq:ClassicalHomotopy} of $h$. Thus $C^\red_2 =
    \hat{C}^\red_2$ in this order of $\nu$. The corrections to the
    terms $\hat{C}^\red_r$ start only in order $r \ge 3$. To compute
    $\hat{C}^\red_2$ we need the second covariant derivatives $\nabla
    \D \pr^*u$ of pull-backs $\pr^*u$. Here we obtain
    \begin{equation}
        \label{eq:NablaXhordpruYhor}
        \left(\nabla_{X^\hor} \D \pr^*u\right) \left(Y^\hor\right)
        = \pr^*\left(\left(\nablaM_X \D u\right) (Y)\right),
    \end{equation}
    \begin{equation}
        \label{eq:NablaXhordpruXH}
        \left(\nabla_{X^\hor} \D \pr^*u\right) \left(X_H\right)
        = - \pr^*\left(\D u (\tau X)\right)
        = \left(\nabla_{X_H} \D \pr^*u\right) \left(X^\hor\right),
    \end{equation}
    \begin{equation}
        \label{eq:NablaXhordpruS}
        \left(\nabla_{X^\hor} \D \pr^*u\right) \left(S\right)
        = - \frac{1}{2} \pr^*\left(\D u(X)\right)
        = \left(\nabla_{S} \D \pr^*u\right) \left(X^\hor\right),
    \end{equation}
    \begin{equation}
        \label{eq:NablaXHdrpuXH}
         \left(\nabla_{X_H} \D \pr^*u\right) \left(X_H\right)
         = \pr^*\left(\D u(\mathsf{V})\right),
     \end{equation}
     \begin{equation}
         \label{eq:NabladpruNull}
         \left(\nabla_{X_H} \D \pr^*u\right) \left(S\right)
         =
         \left(\nabla_{S} \D \pr^*u\right) \left(X_H\right)
         =
         \left(\nabla_{S} \D \pr^*u\right) \left(S\right)
         =
         0.
     \end{equation}
     Inserting this into the general expression \eqref{eq:Czwei} for
     $C_2$ and using $\pr^* \hat{C}^\red_2 (u, v) = H^2 C_2(\pr^*u,
     \pr^*v)$ gives the result \eqref{eq:CTwoRed}. From this, the last
     statement follows directly as the star product $\starred$ is of
     Weyl type and the only second order terms in $C^\red_2$ are
     described by using $\nablaM$, see
     \cite[Prop.~3.1]{gutt.rawnsley:2003a}.
     \qed
\end{proof}

%
%


\end{document}